\documentclass{amsart}
\usepackage{amsfonts,amsmath,amsthm,amssymb}
\usepackage[linesnumbered,ruled,noend]{algorithm2e}
\usepackage{comment}
\usepackage{float}
\usepackage{latexsym}
\usepackage{graphicx}
\usepackage{multicol}
\usepackage{color}
\usepackage{pifont}
\usepackage{enumerate}
\usepackage{lmodern}
\usepackage{makecell}
\usepackage{adjustbox}
\usepackage[backref,pdftex=true,colorlinks,linkcolor=blue,citecolor=red,pdfstartview=FitH]{hyperref}
\usepackage{setspace}
\usepackage[all]{xy}
\usepackage{pb-diagram,pb-xy}
\usepackage{pdflscape}
\usepackage{chngpage}
\usepackage{breakurl}

\usepackage{tabu}

\makeatletter
\def\BState{\State\hskip-\ALG@thistlm}
\makeatother

\theoremstyle{plain}
\newtheorem{theorem}{Theorem}[section]
\newtheorem{lemma}[theorem]{Lemma}

\theoremstyle{definition}

\newtheorem*{example*}{Example}

\DeclareMathOperator{\Rad}{Rad}
\newcommand{\F}{\mathbb{F}}
\newcommand{\Q}{\mathbb{Q}}
\newcommand{\Z}{\mathbb{Z}}

\newcommand{\OO}{\mathcal{O}}

\newcommand{\cP}{\mathcal{P}}
\newcommand{\cE}{\mathcal{E}}

\DeclareMathOperator{\ord}{ord}

\DeclareMathOperator{\Norm}{Norm}

\newcommand{\fq}{\mathfrak{q}}
\newcommand{\sS}{\mathfrak{S}}

\begin{document}

\title[]{On perfect powers that are sums of cubes of a seven term arithmetic progression}

\author{Alejandro Arg\'{a}ez-Garc\'{i}a}
\address{Facultad de Ingenier\'{i}a Qu\'{i}mica, Universidad Aut\'{o}noma de Yucat\'{a}n. Perif\'{e}rico Norte Kil\'{o}metro 33.5, Tablaje Catastral 13615 Chuburna de Hidalgo Inn, M\'{e}rida, Yucat\'{a}n, M\'{e}xico. C.P. 97200 ; 
Escuela Nacional de Estudios Superiores - M\'{e}rida, UNAM,  Calle 7-B No. 227 por 20 y 22-A Colonia Juan B. Sosa M\'{e}rida, Yucat\'{a}n, M\'{e}xico C.P. 97205 }
\email{alejandro.argaez@correo.uady.mx}

\author{Vandita Patel}
\address{School of Mathematics, University of Manchester, Oxford Road, Manchester M13 9PL, United Kingdom
}
\email{vandita.patel@manchester.ac.uk}

\date{\today}


\begin{abstract}
We prove that the equation $(x-3r)^3+(x-2r)^3 + (x-r)^3 + x^3 + (x+r)^3 + (x+2r)^3+(x+3r)^3= y^p$ only has solutions which satisfy $xy=0$ for $1\leq r\leq 10^6$ and $p\geq 5$ prime. This article complements the work on the equations  $(x-r)^3 + x^3 + (x+r)^3 = y^p$ in \cite{ArgaezPatel} and $(x-2r)^3 + (x-r)^3 + x^3 + (x+r)^3 + (x+2r)^3= y^p$  in \cite{Argaez5cubes}. The methodology in this paper makes use of the Primitive Divisor Theorem due to Bilu, Hanrot and Voutier for a complete resolution of the Diophantine equation. 

\end{abstract}
\maketitle

\section{Introduction}

Finding perfect powers that are sums of terms in an arithmetic progression has received much interest lately; recent contributions can be found in \cite{Argaez5cubes}, \cite{ArgaezPatel}, \cite{BGP}, \cite{BPS1}, \cite{BPS2}, \cite{BPSS}, \cite{KoutsianasPatel}, 
\cite{KunduPatel}, \cite{Patelthesis}, \cite{Patelsquares}, \cite{PatelSiksek}, \cite{P}, \cite{P2},
\cite{Soydan},\cite{Van}, \cite{ZB}, \cite{Zhang}, \cite{Zhang2}. 
This paper aims to demonstrate an application of Theorem 1 in \cite{Patel19} to the problem of finding perfect powers in sums of like powers, thus adding to the diverse range of methodologies currently being used to tackle such problems.

\medskip

In this paper, we prove the following: 
\begin{theorem}\label{thm:main}
Let $p \ge 5$ be a prime. The equation
\begin{equation}\label{eq:main}
(x-3r)^3+(x-2r)^3 + (x-r)^3 + x^3 + (x+r)^3 + (x+2r)^3+(x+3r)^3= y^p  
\end{equation}
with $x,r,y,p \in \Z$, $\gcd(x,r)=1$ and $0 < r\le 10^6$ only has solutions which satisfy $xy=0$.
\end{theorem}

The restriction $\gcd(x,r)=1$ is natural one, for otherwise it is easy to construct artificial solutions by scaling. 

\medskip

This paper follows ideas presented in \cite{Patelthesis}, \cite{BPS2}, \cite{ArgaezPatel} and \cite{Argaez5cubes}, to find perfect powers that are sums of terms in an arithmetic progression, whereby the main techniques used to resolve such equations include a result of Mignotte based on linear form in logarithms (\cite[Chapter~12, p.~423]{Cohen2}), the method of Chabauty (\cite{SiksekChab,Stoll, MP}), the theorem due to Bilu, Hanrot and Voutier on primitive divisors (\cite{BHV}), as well as various elementary techniques. 
The novel approach in this paper is that 
we are able to apply a computationally efficient technique
 using Lehmer pairs and sequences due to recent work by the second author \cite{Patel19},
thus leading to a full resolution of \eqref{eq:main} for 
$1\leq r\leq 10^6$ and $p \geq 5$ a prime.
We note here that the approaches taken in \cite{ArgaezPatel, Argaez5cubes} are not sufficient to resolve these cases, as the Thue equations have large coefficients and/or large degree. Similarly, the approach taken in \cite{Patel19} is insufficient by itself. A full resolution of this problem requires a combination of techniques from \cite{ArgaezPatel, Argaez5cubes} and \cite{Patel19}. 

\section{Background}\label{sec:background}

Here, we record some essential theorems and lemmas which are necessary to carry out the computations in section~\ref{sec:tablecomps} in order to prove Theorem~\ref{thm:main}. We emphasise that there is nothing new in this section, and we briefly outline key results and techniques here. Full proofs can be found in \cite{Patelthesis, BPS2, ArgaezPatel}. 

\medskip

We first apply a descent to equation~\eqref{eq:main} in Section~\ref{case12}. 
We are left with equations of the form:
\begin{equation}\label{eqn:form}
aw_2^{p} - bw_1^{2p} = cr^2
\end{equation}
where $p$ is an odd prime and $a,b,c$ are positive integers
satisfying $\gcd(a,b,c)=1$.

Now, we state a theorem due to Mignotte ([\cite{Cohen2}, Chapter 12, p. 423]), which is essential in  providing an upper bound for our prime exponent $p$, thus enabling us to perform a finite computation.

\begin{theorem}\label{thm:mignotte}(Mignotte)
Assume that the exponential Diophantine inequality 
\[
\mid ax^n-by^n\mid\leq c, \qquad \text{with } a,b,c\in\Z_{\geq 0} \text{ and } a\neq b
\]
has a solution in strictly positive integers $x$ and $y$ with $\max\{x,y\}>1$. Let $A=\max\{a,b,3\}$. Then 
\[
n\leq \max\left\lbrace3\log(1.5\mid c/b\mid),\dfrac{7400\log A}{\log\left(1+\log A/\log(\mid a/b\mid) \right)} \right\rbrace.
\]
\end{theorem}

\subsection{Criteria for eliminating equations of signature $(p,2p,2)$} \label{sec:criteria}

Steps to rule out \eqref{eqn:form} having a nontrivial solution for a fixed prime $p\geq 5$ was previously presented in \cite{Patelthesis, BPS2, ArgaezPatel}. We start with the following simple yet effective lemma, inspired from the work of Sophie Germain (see \cite{Patelthesis} for more history and context around the work of Sophie Germain on Fermat's Last Theorem) that provides a computationally efficient criteria to check in order to deduce the nonexistence of solutions. In the majority of cases, this lemma is highly successful in deducing the nonexistence of solutions. 

\begin{lemma}\label{lem:Sophiecriterion}
Let $p \ge 3$ be a prime. Let $a$, $b$ and $c$ be positive integers such that 
$\gcd(a,b,c)=1$. Let $q=2k p+1$ be a prime that does
not divide $a$. Define
\begin{equation}\label{eqn:mu}
\mu(p,q)=\{ \eta^{2p} \; : \; \eta \in \F_q \}
=\{0\} \cup \{ \zeta \in \F_q^* \; : \; \zeta^{k}=1\}
\end{equation}
and
$$
B(p,q)=\left\{ \zeta \in \mu(p,q) \; : \; ((b \zeta+c)/a)^{2k} \in \{0,1\} \right\} \, .
$$
If $B(p,q)=\emptyset$, then equation~\eqref{eqn:form} does not have integral solutions.
\end{lemma}

\subsection{Local Solubility}\label{sub:locsol}

In this section, we outline a classical local solubility method, which when applied, can conclude the nonexistence of solutions for a particular tuple $(a,b,c,p)$ in equation \eqref{eqn:form}.


Recall the condition $\gcd(a,b,c)=1$ in \eqref{eqn:form}. Let $g=\Rad(\gcd(a,c))$ and suppose that $g > 1$.  Then $g \mid w_1$, and we can write $w_1=g w_1^\prime$. Thus
\[
a w_2^p- b g^{2 p} {w_1^\prime}^{2p}=c. 
\]
Removing a factor of $g$ from the coefficients, we obtain
\[
a^\prime {w_2}^p - b^\prime {w_1^\prime}^{2p}=c^\prime,
\]
where $a^\prime=a/c$ and $c^\prime=c/g<c$. Similarly, if $h=\gcd(b,c)>1$, we obtain
\[
a^\prime {w_2^\prime}^p - b^\prime {w_1}^{2p}=c^\prime,
\]
where $c^\prime=c/h<c$. Applying these operations repeatedly, we arrive at an equation of the form
\begin{equation}\label{eqn:predescent}
A \rho^p-B \sigma^{2p}=C
\end{equation}
where $A$, $B$, $C$ are now pairwise coprime. A necessary condition for the existence of solutions is that for any odd prime $q \mid A$, the residue $-BC$ modulo $q$ is a square.
Besides this basic test, we also check for local solubility at the primes dividing $A$, $B$, $C$, and all primes $q \le 19$.

\subsection{Descent}\label{sub:furtherdesc}

If local techniques previously presented fail to rule out solutions 
to equation \eqref{eqn:form} for particular coefficients and exponent $(a,b,c,p)$ 
then we may perform a further descent to rule out solutions. 
With $A$, $B$, $C$ as in \eqref{eqn:predescent}
we let
\[
B^\prime=\prod_{\text{$\ord_q(B)$ is odd}} q.
\]
Thus $B B^\prime=v^2$. Write $A B^\prime=u$ and $C B^\prime=m n^2$
with $m$ squarefree. Rewrite \eqref{eqn:predescent}
as
\[
(v \sigma^p+n \sqrt{-m})(v \sigma^p-n \sqrt{-m})=u \rho^p. 
\]
Let $K=\Q(\sqrt{-m})$ and $\OO$ be its ring of integers. Let $\sS$ contain the prime ideals of $\OO$ that divide $u$ or $2n \sqrt{-m}$. Clearly 
$(v\sigma^p+n \sqrt{-m}) {K^*}^p$ belongs to the ``$p$-Selmer group''
\[
K(\sS,p)=\{\epsilon \in K^*/{K^*}^p \; : \; 
\text{$\ord_\cP(\epsilon) \equiv 0 \mod{p}$ for all $\cP \notin \sS$}
\}.
\]
This is an $\F_p$-vector space of finite dimension can be computed by \texttt{Magma} using the command \texttt{pSelmerGroup}. Let
\[
\cE=\{ \epsilon \in K(\sS,p) \; : \; \Norm(\epsilon)/u \in {\Q^*}^p \}.
\]
It follows that
\begin{equation}\label{eqn:furtherdescent}
v \sigma^p+n \sqrt{-m}=\epsilon \eta^p,
\end{equation}
where $\eta \in K^*$ and $\epsilon \in \cE$.


\begin{lemma}\label{lem:valuative}
Let $\fq$ be a prime ideal of $K$. Suppose one of the following
holds:
\begin{enumerate}[$(i)$]
\item $\ord_\fq(v)$, $\ord_\fq(n\sqrt{-m})$, $\ord_\fq(\epsilon)$
are pairwise distinct modulo $p$;
\item $\ord_\fq(2v)$, 
$\ord_\fq(\epsilon)$, $\ord_\fq(\overline{\epsilon})$
are pairwise distinct modulo $p$;
\item $\ord_\fq(2 n \sqrt{-m})$, 
$\ord_\fq(\epsilon)$, $\ord_\fq(\overline{\epsilon})$
are pairwise distinct modulo $p$.
\end{enumerate}
Then there is no $\sigma \in \Z$ and $\eta \in K$ satisfying \eqref{eqn:furtherdescent}.
\end{lemma}

\begin{lemma}\label{lem:furtherdescent}
Let $q=2k p+1$ be a prime. 
Suppose $q\OO=\fq_1 \fq_2$ where $\fq_1$, $\fq_2$
are distinct, and such that $\ord_{\fq_j}(\epsilon)=0$
for $j=1$, $2$. Let 
\[
\chi(p,q)=\{ \eta^p \; : \; \eta \in \F_q \}.
\]
Let
\[
C(p,q)=\{\zeta \in \chi(p,q) \; : \;
((v \zeta+n\sqrt{-m})/\epsilon)^{2k} \equiv \text{$0$ or $1 \mod{\fq_j}$
for $j=1$, $2$}\}.
\]
Suppose $C(p,q)=\emptyset$. 
Then there is no $\sigma \in \Z$
and $\eta \in K$ satisfying \eqref{eqn:furtherdescent}.
\end{lemma}

\subsection{Thue equations}
For the handful of remaining equations where we are unable to deduce nonexistence of solutions, we let $\sigma = w_2$ and $\tau = w_1^2$ in \eqref{eqn:form} to get:
\[
a\sigma^p - b \tau^p = cr^2
\]
For a fixed value of $r$, we note that this is a \emph{Thue equation}. We use \texttt{Magma's} Thue solver \cite{magma} and \texttt{PARI/GP's} \textit{thueinit}, \textit{thue} commands \cite{PARI2, sagemath}, as the final test to determine whether the equation has solutions.

\section{The initial descent and twelve cases }\label{case12}
We rewrite equation~\eqref{eq:main} as $7x(x^2+12r^2)= y^p$. Since $7 \mid y$, we let $y = 7w$ to obtain:
\[
x(x^2+12r^2) = 7^{p-1}w^p.
\]
We note that $\gcd(x, x^2 + 12r^2) \in \{1,2,3,4,6,12\}$ depending on whether $2,3,4,6$ and $12$ divides $x$ or not. 
This leads us to consider twelve cases. We apply a simple descent argument in each case and the results are summarised in the following table.

\medskip 

\begin{adjustbox}{width =\textwidth}
\begin{tabular}{|c|c|c|c|}
\hline
{\bf Case} & {\bf Conditions on $x$} & {\bf Descent equations} & {\bf Ternary Equation}\\
\hline\hline
$1$ &  $7 \nmid x$ and $12 \nmid x$ &
\begin{tabular}{@{}c@{}}$x=w_1^p$ \\ $x^2+12r^2=7^{p-1}w_2^p$ \end{tabular} & $7^{p-1}w_2^p-w_1^{2p}=12r^2$ \\

\hline
$2$ &  $7 \nmid x$ and $2 \nmid x$, $3\mid x$ &
\begin{tabular}{@{}c@{}}$x=3^{p-1}w_1^p$ \\ $x^2+12r^2=3\cdot7^{p-1}w_2^p$ \end{tabular} & $7^{p-1}w_2^p-3^{2p-3}w_1^{2p}=4r^2$ \\

\hline
$3$ &  $7 \nmid x$ and $4 \mid x$, $3\nmid x$ &
\begin{tabular}{@{}c@{}}$x=2^{p-2}w_1^p$ \\ $x^2+12r^2=4\cdot7^{p-1}w_2^p$ \end{tabular} & $7^{p-1}w_2^p-2^{2p-6}w_1^{2p}=3r^2$ \\

\hline
$4$ &  $7 \nmid x$ and $12 \mid x$ &
\begin{tabular}{@{}c@{}}$x=2^{p-2}\cdot3^{p-1}w_1^p$ \\ $x^2+12r^2=12\cdot7^{p-1}w_2^p$ \end{tabular} & $7^{p-1}w_2^p-2^{2p-6}\cdot3^{2p-3}w_1^{2p}=r^2$ \\

\hline 
$5$ &  $7 \nmid x$ and $2 \mid x$, $3,4\nmid x$ &
\begin{tabular}{@{}c@{}}$x=2^{p-1}w_1^p$ \\ $x^2+12r^2=2\cdot7^{p-1}w_2^p$ \end{tabular} & $7^{p-1}w_2^p-2^{2p-3}w_1^{2p}=6r^2$ \\

\hline 
$6$ &  $7 \nmid x$ and $6 \mid x$, $4\nmid x$ &
\begin{tabular}{@{}c@{}}$x=6^{p-1}w_1^p$ \\ $x^2+12r^2=6\cdot7^{p-1}w_2^p$ \end{tabular} & $7^{p-1}w_2^p-6^{2p-3}w_1^{2p}=2r^2$ \\

\hline
$7$ &  $7 \mid x$ and $12 \nmid x$ &
\begin{tabular}{@{}c@{}}
$x=7^{p-1}w_1^p$ \\  $x^2+12r^2=w_2^p$ 
\end{tabular} 
& $w_2^p-7^{2p-2}w_1^{2p}=12r^2$\\

\hline 
$8$ &  $7 \mid x$ and $2 \nmid x$, $3\mid x$ &
\begin{tabular}{@{}c@{}}$x=3^{p-1}\cdot7^{p-1}w_1^p$ \\ $x^2+12r^2=3w_2^p$ \end{tabular} & $w_2^p-3^{2p-3}\cdot7^{2p-2}w_1^{2p}=4r^2$ \\

\hline 
$9$ &  $7 \mid x$ and $4 \mid x$, $3\nmid x$ &
\begin{tabular}{@{}c@{}}$x=2^{p-2}\cdot7^{p-1}w_1^p$ \\ $x^2+12r^2=4w_2^p$ \end{tabular} & $w_2^p-2^{2p-6}\cdot7^{2p-2}w_1^{2p}=3r^2$ \\

\hline 
$10$ &  $7 \mid x$ and $12 \mid x$ &
\begin{tabular}{@{}c@{}}$x=2^{p-2}\cdot3^{p-1}\cdot7^{p-1}w_1^p$ \\ $x^2+12r^2=12w_2^p$ \end{tabular} & $w_2^p-2^{2p-6}\cdot3^{2p-3}\cdot7^{2p-2}w_1^{2p}=r^2$ \\

\hline 
$11$ &  $7 \mid x$ and $2 \mid x$, $3,4\nmid x$ &
\begin{tabular}{@{}c@{}}$x=2^{p-1}\cdot7^{p-1}w_1^p$ \\ $x^2+12r^2=2w_2^p$ \end{tabular} & $w_2^p-2^{2p-3}\cdot7^{2p-2}w_1^{2p}=6r^2$ \\

\hline
$12$ &  $7 \mid x$ and $6 \mid x$, $4\nmid x$ &
\begin{tabular}{@{}c@{}}$x=6^{p-1}\cdot7^{p-1}w_1^p$ \\ $x^2+12r^2=6w_2^p$ \end{tabular} & $w_2^p-6^{2p-3}\cdot7^{2p-2}w_1^{2p}=2r^2$ \\

\hline
\end{tabular}
\label{table:Cases}
\end{adjustbox}

\medskip

We note that cases $5$, $6$, $11$ and $12$  all lead to an immediate contradiction; by checking the valuation of $2$ on both sides of the ternary equation and noting the fact that $p \geq 5$, we contradict our assumption of $\gcd(x,r)=1$.


\section{Solving cases $1$--$4$}\label{sec:tablecomps}
In this section, we apply the Theorems and Lemmas of Section~\ref{sec:background} which fully resolve descent cases $1$--$4$ of equation~\eqref{eq:main} in order to prove Theorem~\ref{thm:main}.

\subsection{Case 1.}

For $p=5$ we have $7^4w_2^5-w_1^{10}=3(2r)^2$. Letting $X=w_2/w_1^2$ and $Y=6r/w_1^5$, we obtain the hyperelliptic curve 
\[
Y^2=3\cdot7^4X^5-3
\]
whose Jacobian has rank 1. The Chabauty implementation give us $C(\Q)=\{\infty\}$.

By Theorem~$\ref{thm:mignotte}$, for $|r|\leq 4.9\times 10^{1502}$, we can bound $p\leq 20775$. Thus, for $7\leq p \leq 20775$ and $1\leq r\leq 10^6$ we have the following table.

\[
\begin{adjustbox}{width =\textwidth}
\begin{tabular}{|c|c|c|c|c|}
\hline
\text{Exponent }p & \makecell{Number of eqns\\ surviving \\Lemma \ref{lem:Sophiecriterion}}& \makecell{Number of eqns \\surviving local \\solubility tests} & \makecell{Number of eqns \\surviving \\further descent} & \makecell{Thue eqns \\ not solved  \\ by Magma} \\
\hline
7 & 37679 & 4077 & 3 & 0\\
\hline
11 & 9930 & 5375 & 0 & 0\\
\hline
13 & 3298 & 1405 & 0 & 0\\
\hline
17 & 461 & 253 & 0 & 0\\
\hline
19 & 1507 & 936 & 0 & 0\\
\hline
23 & 13 & 3 & 0 & 0\\
\hline
29 & 8 & 5 & 0 & 0\\
\hline
31 & 29 & 21 & 0 & 0\\
\hline
$37\leq p \leq 20775$& 0 & 0 & 0 & 0\\
\hline
\end{tabular}
\end{adjustbox}
\]

\subsection{Case 2.}

For $p=5$ we have $7^4w_2^5-3^7w_1^{10}=(2r)^2$. Letting $X=w_2/w_1^2$ and $Y=2r/w_1^5$ we obtain the hyperelliptic curve 
\[
Y^2=7^4X^5-3^7
\]
whose Jacobian has rank 2, and we are unable to use Chabauty techniques.

By Theorem~$\ref{thm:mignotte}$, for $|r|\leq 1.9\times 10^{1427}$, we can bound $p\leq 19734$. Thus, for $5\leq p \leq 19734$ and $1\leq r\leq 10^6$ we have the following table.

\[
\begin{adjustbox}{width =\textwidth}
\begin{tabular}{|c|c|c|c|c|}
\hline
\text{Exponent }p & \makecell{Number of eqns\\ surviving \\Lemma \ref{lem:Sophiecriterion}}& \makecell{Number of eqns \\surviving local \\solubility tests} & \makecell{Number of eqns \\surviving \\further descent} & \makecell{Thue eqns \\ not solved  \\ by Magma} \\
\hline
5 & 102681 & 38771 & 819 & 819 \\
\hline
7 & 24526 & 2400 & 0 & 0\\
\hline
11 & 9159 & 4629 & 0 & 0\\
\hline
13 & 3439 & 1804 & 0 & 0\\
\hline
17 & 80 & 51 & 0 & 0\\
\hline
19 & 3136 & 1012 & 0 & 0\\
\hline
23 & 11 & 3 & 0 & 0\\
\hline
29 & 1 & 0 & 0 & 0\\
\hline
31 & 4 & 2 & 0 & 0\\
\hline
$37\leq p \leq 19734$ & 0 & 0 & 0 & 0\\
\hline
\end{tabular}
\end{adjustbox}
\]

Using \texttt{PARI/GP}, we resolved the $819$ equations using the commands \textit{thueinit} and \textit{thue}. No integer solutions were found.

\subsection{Case 3.}

For $p=5$ we have $7^6w_2^5-2^4w_1^{10}=3r^2$. Letting $X=w_2/w_1^2$ and $Y=3r/w_1^5$ we obtain the hyperelliptic curve 
\[
Y^2=3\cdot7^6X^5-2^4\cdot 3
\]
whose Jacobian has rank $0$. Applying Chabauty gives $C(\Q)=\{\infty\}$.

By Theorem~\ref{thm:mignotte}, for $|r|\leq 1.5\times 10^{2105}$, we can bound $p\leq 29101$. Thus, for $7\leq p\leq 29101$ and $1\leq r\leq 10^6$ we have the following table.

\[
\begin{adjustbox}{width =\textwidth}
\begin{tabular}{|c|c|c|c|c|}
\hline
\text{Exponent }p & \makecell{Number of eqns\\ surviving \\Lemma \ref{lem:Sophiecriterion}}& \makecell{Number of eqns \\surviving local \\solubility tests} & \makecell{Number of eqns \\surviving \\further descent} & \makecell{Thue eqns \\ not solved  \\ by Magma} \\
\hline
7  & 29213 & 2969 & 0 & 0\\
\hline
11 & 6484 & 2332 & 0 & 0\\
\hline
13 & 1715 & 786 & 0 & 0\\
\hline
17 & 369 & 206 & 0 & 0\\
\hline
19 & 538 & 262 & 0 & 0\\
\hline
23 & 1 & 0 & 0 & 0\\
\hline
29 & 2 & 1 & 0 & 0\\
\hline 
31 & 5 & 4 & 0 & 0\\
\hline
37$\leq p \leq $29101 & 0 & 0 & 0 & 0 \\
\hline
\end{tabular}
\end{adjustbox}
\]

\subsection{Case 4.}

For $p=5$ we have $7^6w_2^5-2^4\cdot 3^7w_1^{10}=r^2$. Letting $X=w_2/w_1^2$ and $Y=r/w_1^5$ we obtain the hyperelliptic curve
\[
Y^2=7^6X^5-2^4\cdot 3^7
\]
whose Jacobian has rank $1$. The Chabauty implementation gives $C(\Q)=\{\infty\}$.

By Theorem~\ref{thm:mignotte}, for $|r|\leq 1.37\times 10^{4664}$, we can bound $p\leq 64461$. Hence, for $7\leq p\leq 64461$ and $1\leq r\leq 10^6$ we have the following table.

\[
\begin{adjustbox}{width =\textwidth}
\begin{tabular}{|c|c|c|c|c|}
\hline
\text{Exponent }p & \makecell{Number of eqns\\ surviving \\Lemma \ref{lem:Sophiecriterion}}& \makecell{Number of eqns \\surviving local \\solubility tests} & \makecell{Number of eqns \\surviving \\further descent} & \makecell{Thue eqns \\ not solved  \\ by Magma} \\
\hline
7  & 18908 & 1940 & 0 & 0\\
\hline
11 & 1384  & 434  & 0 & 0\\
\hline
13 & 479   & 177  & 0 & 0\\ 
\hline
17 & 366   & 173  & 0 & 0\\
\hline
19 & 365   & 184  & 0 & 0\\
\hline
23 & 5     & 1    & 0 & 0\\
\hline
29 & 3     & 3    & 0 & 0\\
\hline
31 & 14    & 9    & 0 & 0\\
\hline
37 & 1     & 0    & 0 & 0\\
\hline
$41\leq p \leq 64461$ & 0 & 0 & 0 & 0\\
\hline
\end{tabular}
\end{adjustbox}
\]

%

%

%

\section{Solving cases $7$--$10$}\label{sec:7to12}
For cases $7$--$10$, we rely on techniques developed in \cite{Patel19}.
Using theorems and lemmas of Section~\ref{sec:background} lead to resolving Thue equations with extremely large coefficients and/or high degree. In order to make computations tractable, we now outline an alternative approach to deal with cases $7$--$10$.

\subsection{Primitive prime divisors of Lehmer sequences}\label{sec:Lehmer}

\begin{theorem}\label{thm:main2}
Let $C_1 \ge 1$ be a squarefree integer and $C_2$ a positive integer. Write 
$C_1C_2 = cd^2$ where $c$ is squarefree. We assume that $C_1C_2 \not\equiv 7
\pmod{8}$.  Let $p$ be a prime for which 
\begin{equation}\label{eqn:main2}
C_1x^2 + C_2 = y^p, \quad x,~y \in \Z^{+}, \quad \gcd(C_1x^2,C_2,y^p)=1,
\end{equation}
has a solution $(x,y)$.
Then either, 
\begin{enumerate}[(i)]
\item $p \le 13$, or
\item $p$ divides the class number of $\Q(\sqrt{-c})$, or
\item $p \mid \left( q - \left(\frac{-c}{q}\right)\right)$, where $q$ is some
prime $q\mid d$ and $q\nmid 2c$.  
\end{enumerate}
\end{theorem}

We follow \textbf{Case I} in Section~$6$ of \cite{Patel19} which outlines effective methods to solve equation~\eqref{eqn:main2} for fixed values of $C_1, C_2, p$. 
Theorem~\ref{thm:main2} is crucial to gain good bounds for $p$. For the convenience of the reader, we shall reproduce the key computational steps used to prove Theorem~\ref{thm:main} here. For the full exposition and intricate details of computations, see \cite{Patel19}.

\medskip

Let $C_1$, $C_2$ be positive integers, with $C_1$ squarefree. Let $\gcd(C_1, C_2)=1$ and suppose that $C_1C_2 \not \equiv 7 \pmod{8}$. We write $C_1C_2 = cd^2$ where $c, d$ are positive integers and $c$ is squarefree.
Theorem~\ref{thm:main2} gives a list of possible odd prime
exponents $p$ for which \eqref{eqn:main2} might have solutions.

Let $K=\Q(\sqrt{-3})$ be a number field with $\OO_K$ it's ring of integers.
We note here that $c=3$ in the notation of Theorem~\ref{thm:main2}. Furthermore, we note that the class number of $K$ is $1$. Let $p\geq 5$ be a prime.
Let $\gamma \in \OO_K$ such that $\gamma = a + b\sqrt{-3}$ for some integers $a,b$,
and let $\bar\gamma$ be it's conjugate.

Fix a value of $b$ dividing $d$.  
To determine the solutions we merely have to determine the possible
values of $a$ corresponding to each $b \mid d$.
We write an explicit polynomial $g_b \in \Z[X]$
whose integer roots contain all the possible values of $a$
corresponding to $b$.

Fix $s \mid d$. 
Since $-3 \not \equiv 1 \pmod{4}$, we let
\begin{equation}\label{eqn:roots}
g_b(X)=\frac{ (X+b\sqrt{-3})^p - (X-b\sqrt{-3})^p}{2b \sqrt{-3}}
\; - \; \frac{d \cdot C_1^{(p-1)/2}}{b}.
\end{equation}
Clearly $g_b \in \Z[X]$. 
Moreover,
\[
g_b(a)  =
\frac{ \gamma^p - \overline{\gamma}^p}{\gamma-\overline{\gamma}}
\; - \; \frac{d \cdot C_1^{(p-1)/2}}{b}
=0.
\]
Finding the integer roots $g_b(X)$ gives associated values of $a$ to the specified $b$, hence $\gamma$ is known. 
Since $y = \gamma\bar\gamma/C_1$, hence $y= (a^2 + 3b^2)/C_1$, the possible values for $y$ are now easily known. Using equation~\eqref{eqn:main2}, we finally determine possible values of $x$, hence all possible solutions to equation~\eqref{eqn:main2}.

\subsection{Case 7.}
One of the descent equations is:
\begin{equation}\label{eqn:case7}
x^2 + 12r^2 = w_2^p
\end{equation}

We let $C_1 = 1$, $C_2 = 12r^2$ and we see that $C_1C_2 \equiv 0,3,4 \pmod{8} \not\equiv 7 \pmod{8}$. Thus if \eqref{eqn:case7} has a solution, then we have $5 \leq p \leq 13$ or  $p \mid \left( q - \left(\frac{-3}{q}\right)\right)$, where $q$ is some prime $q\mid r$ and $q\nmid 6$. For fixed values of $1\leq r \leq 10^6$ (hence fixed $C_1$, $C_2$, $p$), finding the roots of the polynomial \eqref{eqn:roots}, leads to solutions $(x,w_2,p)$. Thus to obtain solutions $(x,y,p)$ to the original equation, we simply note that 
$y=7w_1w_2$ where $w_1$ is easily deduced from the descent equation $x=7^{p-1}w_1^p$ in Section~\ref{case12}. The computation yielded no such solutions. 


\subsection{Case 8.}
One of the descent equations is:
\[
x^2 + 12r^2 = 3w_2^p.
\]
We let $x = 3X$ and obtain
\begin{equation}\label{eqn:case8}
3X^2 + 4r^2 = w_2^p.
\end{equation}

We let $C_1 = 3$, $C_2 = 4r^2$ and we see that $C_1C_2 \equiv 0,3,4 \pmod{8} \not\equiv 7 \pmod{8}$. Thus if \eqref{eqn:case8} has a solution, then we have $5 \leq p \leq 13$ or  $p \mid \left( q - \left(\frac{-3}{q}\right)\right)$, where $q$ is some prime $q\mid r$ and $q\nmid 6$. For fixed values of $1\leq r \leq 10^6$ (hence fixed $C_1$, $C_2$, $p$), finding the roots of the polynomial \eqref{eqn:roots}, leads to solutions $(X,w_2,p)$. Thus to obtain solutions $(x,y,p)$ to the original equation, we simply note that $x = 3X$ and  
$y=7w_1w_2$ where $w_1$ is easily deduced from the descent equation $x=3^{p-1}7^{p-1}w_1^p$ in Section~\ref{case12}. The computation yielded no such solutions. 

\subsection{Case 9.}
One of the descent equations is:
\[
x^2 + 12r^2 = 4w_2^p.
\]
We let $x = 2X$ and obtain
\begin{equation}\label{eqn:case9}
X^2 + 3r^2 = w_2^p.
\end{equation}

We let $C_1 = 1$, $C_2 = 3r^2$ and we see that $C_1C_2 \equiv 0,3,4 \pmod{8} \not\equiv 7 \pmod{8}$. Thus if \eqref{eqn:case9} has a solution, then we have $5 \leq p \leq 13$ or  $p \mid \left( q - \left(\frac{-3}{q}\right)\right)$, where $q$ is some prime $q\mid r$ and $q\nmid 6$. For fixed values of $1\leq r \leq 10^6$ (hence fixed $C_1$, $C_2$, $p$), finding the roots of the polynomial \eqref{eqn:roots}, leads to solutions $(X,w_2,p)$. Thus to obtain solutions $(x,y,p)$ to the original equation, we simply note that $x = 2X$ and  
$y=7w_1w_2$ where $w_1$ is easily deduced from the descent equation $x=2^{p-2}7^{p-1}w_1^p$ in Section~\ref{case12}. The computation yielded no such solutions.

\subsection{Case 10.}
One of the descent equations is:
\[
x^2 + 12r^2 = 12w_2^p.
\]
We let $x = 6X$ and obtain
\begin{equation}\label{eqn:case10}
3X^2 + r^2 = w_2^p.
\end{equation}

We let $C_1 = 3$, $C_2 = r^2$ and we see that $C_1C_2 \equiv 0,3,4 \pmod{8} \not\equiv 7 \pmod{8}$. Thus if \eqref{eqn:case10} has a solution, then we have $5 \leq p \leq 13$ or  $p \mid \left( q - \left(\frac{-3}{q}\right)\right)$, where $q$ is some prime $q\mid r$ and $q\nmid 6$. For fixed values of $1\leq r \leq 10^6$ (hence fixed $C_1$, $C_2$, $p$), finding the roots of the polynomial \eqref{eqn:roots}, leads to solutions $(X,w_2,p)$. Thus to obtain solutions $(x,y,p)$ to the original equation, we simply note that $x = 6X$ and  
$y=7w_1w_2$ where $w_1$ is easily deduced from the descent equation $x=2^{p-2}3^{p-1}7^{p-1}w_1^p$ in Section~\ref{case12}. The computation yielded no such solutions.


\bibliographystyle{plain}

\begin{thebibliography}{}

\bibitem{Argaez5cubes}
A.  Arg\'{a}ez-Garc\'{i}a,
\emph{On perfect powers that are sums of cubes of a five term arithmetic progression},
J. Number Theory, \textbf{201} (2019), 460--472.

\bibitem{ArgaezPatel}
A. Arg\'{a}ez-Garc\'{i}a and V. Patel,
\emph{Perfect powers that are sums of cubes of a three term arithmetic progression},
Journal of Combinatorics and Number Theory (2019) (to be announce).

\bibitem{BGP}
M. A. Bennett, K. Gy\H{o}ry and \'{A}. Pint\'{e}r,
\emph{On the Diophantine equation $1^k+2^k+\cdots+x^k=y^n$},
Compos. Math. \textbf{140} (2004), no. 6, 1417--1431. 

\bibitem{BPS1}
M.\ A.\ Bennett, V.\ Patel and S.\ Siksek,
\emph{Superelliptic equations arising from sums of consecutive powers},
Acta Arith., \textbf{172} (2016), no. 4, 377--393.

\bibitem{BPS2}
M.\ A.\ Bennett, V.\ Patel and S.\ Siksek,
\emph{Perfect powers that are sums of consecutive cubes},
Mathematika, \textbf{63} (2017), no. 1, 230--249.

\bibitem{BPSS} {\sc A. B\'{e}rczes, I. Pink, G. Sava\c{s}, G. Soydan}, On the Diophantine equation $(x+1)^k+ (x+2)^k+ ... + (2x)^k=y^{n},$ {\em J. Number Theory } {\bf 183}(2018), 326-351.

\bibitem{BHV}
Yu.\ Bilu, G.\ Hanrot, and P.\ M.\ Voutier, 
\emph{Existence of primitive divisors of Lucas and Lehmer numbers},
 J. Reine Angew. Math. \textbf{539} (2001), 75--122.

\bibitem{magma}
W.\ Bosma, J.\ Cannon and C.\ Playoust,
\emph{The Magma Algebra System I: The User Language},
J.\ Symb.\ Comp. \textbf{24} (1997), 235--265. (See also \url{http://magma.maths.usyd.edu.au/magma/})



\bibitem{Cohen2}
H.\ Cohen,
\emph{Number theory. {V}ol. {II}. {A}nalytic and modern tools}, Graduate Texts in Mathematics \textbf{240}, Springer, New York, 2007.


\bibitem{KoutsianasPatel}
A. Koutsianas, V. Patel,
\emph{Perfect powers that are sums of squares in a three term arithmetic progression},
Int. J. Number Theory, \textbf{14}(10) (2018).

\bibitem{KunduPatel}
D. Kundu and V. Patel,
\emph{Perfect powers that are sums of squares of an arithmetic progression},
\href{https://arxiv.org/abs/1809.09167}{arXiv:1809.09167}, (2018).


\bibitem{MP} 
W. McCallum and B.\ Poonen, 
\emph{The method of Chabauty and Coleman}, pages 99--117 of \cite{DioBook}.

 \bibitem{PARI2}
    The PARI~Group, PARI/GP version {\tt 2.11.0}, Univ. Bordeaux, 2018,
    \url{http://pari.math.u-bordeaux.fr/}.

\bibitem{Patelthesis}
V.\ Patel,
\emph{Perfect powers that are sums of consecutive like
powers},
Doctoral thesis (2017).


\bibitem{Patelsquares}
V. Patel,
\emph{Perfect powers that are sums of consecutive squares},
C. R. Math. Rep. Acad. Sci. Can., \textbf{40}(2) (2018),  33--38.

\bibitem{Patel19}
V. Patel,
\emph{A Lucas-Lehmer approach to generalised Lebesgue-Ramanujan-Nagell equations},
\href{https://arxiv.org/abs/1910.07453}{arXiv:1910.07453}, (2019).


\bibitem{PatelSiksek}
V.\ Patel and S.\ Siksek,
\emph{On powers that are sums of consecutive like powers},
Research in Number Theory \textbf{3} (2017), 2:7.

\bibitem{P}
\'A. Pint\'er,
\emph{A note on the equation $1^k+2^k+ \cdots + (x-1)^k=y^m$},
Indag. Math. N.S. \textbf{8} (1997), 119--123.

\bibitem{P2}
\'A. Pint\'er,
\emph{On the power values of power sums},
J. Number Theory  \textbf{125} (2007), 412--423.



\bibitem{SiksekChab}
\emph{Chabauty and the Mordell-Weil sieve}, pages 194--224 of
\textit{Advances on superelliptic curves and their applications},  
NATO Sci. Peace Secur. Ser. D Inf. Commun. Secur. \textbf{41}, IOS, Amsterdam, 2015. 

\bibitem{Soydan} {\sc G. Soydan}, On the Diophantine equation   $(x+1)^k+ (x+2)^k+ ... + (lx)^k=y^{n}$, {\em Publ. Math. Debrecen } {\bf 91} (2017), 369-382.

\bibitem{Stoll}
M.\ Stoll, 
\emph{Implementing 2-descent for Jacobians of hyperelliptic curves},
 Acta Arith. \textbf{98} (2001), no.\ 3, 245--277.


\bibitem{Van}
J.\ van Lagen,
\emph{On the sum of fourth powers in arithmetic progression},
\href{https://arxiv.org/abs/1907.12351}{arXiv:1907.12351}, (2019).

\bibitem{ZB}
Z.\ Zhang and M. Bai,
\emph{On the Diophantine equation $(x+1)^2+(x+2)^2 + \cdots +(x+d)^2=y^n$},
Funct. Approx. Comment. Math.  \textbf{49} (2013), 73--77.
\bibitem{sagemath}
Developers, The~Sage,
\emph{{S}agemath, the {S}age {M}athematics {S}oftware {S}ystem ({V}ersion 8.5), 2019, \tt https://www.sagemath.org}

\bibitem{Zhang}
Z.\ Zhang,
\emph{On the Diophantine equation $(x-1)^k+x^k+(x+1)^k=y^n$},
Publ.\ Math.\ Debrecen \textbf{85} (2014), 93--100.

\bibitem{Zhang2}
Z.\ Zhang,
\emph{On the Diophantine equation $(x-d)^4+x^4+(x+d)^4=y^n$},
Int.\ J.\ Number Theory \textbf{} \emph{published online} (2017)., 93--100.






\end{thebibliography}


\end{document}